\magnification=\magstep1
\input amstex
\documentstyle{amsppt}
\voffset=-3pc
\loadbold
\loadmsbm
\loadeufm
\UseAMSsymbols
\baselineskip=12pt
\parskip=6pt

\def\bC{\Bbb C}

\def\bR{\Bbb R}

\def\oP{\overline P}
\def\bP{\Bbb P}
\def\id{\text{id}}

\NoBlackBoxes
\topmatter
\title A new look at adapted complex structures\endtitle
\author L\'aszl\'o Lempert and 
R\'obert Sz\H{o}ke$^{*)}$\footnote""{$^*$Research supported by NSF grant DMS0700281 and OTKA grants T047102, K72537.
The second author acknowledges the hospitality of the 
Purdue Mathematics Department, where most of this research
was done during his visit in 2005.\hfill\break}\endauthor
\rightheadtext{L\'aszl\'o Lempert and R\'obert Sz\H{o}ke }
\leftheadtext{A New Look at Adapted Complex Structures}
\address
Dept.~of Mathematics, Purdue University, West Lafayette, IN\ 47907, USA\endaddress
\address Dept.~of Analysis, Institute of Mathematics,
E\"otv\"os University, P\'azm\'any P. 
s\'et\'any 1/c, Budapest 1117, Hungary\endaddress
\subjclassyear{2000}
\subjclass 32Q15, 53C20,22\endsubjclass
\abstract
Given a closed real analytic Riemannian manifold, we construct 
and study a one parameter family of adapted complex structures on the manifold of its geodesics.
\endabstract
\dedicatory R. Sz. dedicates this paper to his son, B\'alint
\enddedicatory
\endtopmatter
\document

\TagsOnRight
\subhead 1.\ Introduction\endsubhead
Adapted complex structures, also called Grauert tubes, are K\"ahler structures on the tangent or cotangent bundle of a Riemannian manifold $M$.
They first appeared in \cite{GS} and \cite{LSz1}.
Let $T^r M\subset TM$ denote the set of tangent vectors of length $<r$.
By the above papers and by \cite{L, Sz1} one knows that for a closed $M$ an adapted complex structure exists on $T^r M$ for some $r>0$ if and only if the metric of $M$ is analytic; moreover, when it exists, the adapted complex structure is unique.
These results have lately been extended to Koszul manifolds, i.e., manifolds with a connection, in \cite{B, Sz3}.

The thesis of this note is that one can talk about adapted complex and related structures in generality much greater than that of Riemannian or Koszul manifolds, indeed whenever a Lie semigroup acts on a manifold.
Of course, whether these more general structures exist, are unique, or have useful properties, depends on the action.
Our focus here will still be on Riemannian manifolds.
For them the new definition is not quite equivalent to the one in \cite{GS, LSz1}.
Rather, there is a family of equally natural K\"ahler structures on (subdomains of) $TM$, parametrized by $s\in\bC\backslash\bR$.
The K\"ahler manifolds thus obtained constitute the fibers of a holomorphic fibration over $\bC\backslash\bR$, and the adapted complex structure of \cite{GS, LSz1} corresponds to the fiber 
over $s=i$.---It is possible to extend the fibration to a 
fibration over $\bC$; however, the fibers over $\bR$ will be real polarized rather than K\"ahler.
Thus one is led to the notion of adapted polarizations, of which an adapted complex structure is just an extreme example.

The idea that the adapted complex structures of \cite{GS, LSz1}
give rise to other K\"ahler structures as well is as old as the whole subject.
More to the point, the papers [FMMN1--2] consider a one (real) parameter family of K\"ahler structures on the cotangent bundle of a compact Lie group, that degenerates to a real polarization; this family is then used to explain geometrically the so called Bargmann--Segal--Hall transformation of \cite{H1-2}.
The papers themselves make no explicit connection with adapted complex structures, but the family considered there is the restriction of our family of adapted polarizations to the positive imaginary axis.
Rather recently \cite{HK} pointed out that for a general closed real analytic Riemannian manifold the original adapted complex structure is the analytic continuation to $i$ of a real family of real polarizations.
Undoubtedly, the authors were aware that continuation to other values $s\in\bC\backslash\bR$ also yields K\"ahler structures.

The novelty of the present note is first that all those K\"ahler structures and real polarizations can be derived from one principle; second that these structures, taken together, constitute a fiber bundle.
The results obtained will be used in the companion paper [LSz2] to study the uniqueness of geometric quantization.

\subhead 2.\ Polarized manifolds\endsubhead
A polarization of a smooth manifold $N$ is given by a smooth, involutive, complex subbundle $P\subset\bC\otimes TN$, of rank $m=(1/2)\dim N$.
Involutivity means that the bracket of sections of $P$ is again a section of $P$.
This definition is more general than the one, say, in [W], but in our context this is the natural one.
Sometimes even more general structures have to be considered, where the rank condition is omitted; these are the involutive manifolds.

A polarization is real if $\oP=P$; it is equivalent to the datum of an $m$ dimensional foliation of $N$.
A polarization is complex if $P$ and $\oP$ are transverse; this 
one is equivalent to a complex structure on $N$.
In the former case $P$ consists of tangents to the leaves, in the latter $P$ is the bundle of $(1,0)$ vectors.
A smooth map $f$ of polarized manifolds $(M,Q)\to (N,P)$ is called polarized if $f_*Q\subset P$.

Consider now a smooth manifold $N$ on which a Lie semigroup $G$ acts smoothly on the right.
Fix a polarization of $G$.

\definition{Definition 1}A polarization of $N$ is called adapted (to the polarization of $G$) if for every $x\in N$ the map $G\ni g\mapsto xg\in N$ is polarized.
\enddefinition

All of the above makes sense and works the same
when $N$ is a manifold with boundary, except that the connection
between real polarizations and foliations becomes tenuous (for
there is more than one way to define foliations on such $N$).

\subhead 3.\ The affine semigroup\endsubhead
This is the Lie semigroup $G$ of affine transformations $\bR\to\bR$, $t\mapsto gt=a+bt$.
Here $a=a(g), b=b(g)$ serve as global coordinates on $G$, and identify it with $\bR^2$.
Denote by $L_g,R_g$ left and right translations of $G$.
By a left invariant polarization $Q$ of $G$ we mean one for which $L_g\colon G\to G$ is polarized for every $g\in G$.
In the identification $G\cong\bR^2$ invariance means that the fibers $Q_g$, $g\in G$, are Euclidean translates of each other.
Hence associating with a polarization $Q$ of $G$ the complex line $Q_{\text{id}}\subset\bC\otimes T_{\id} G$ yields a bijection between the set of left invariant polarizations and $\bP(\bC\otimes T_{\id} G)\approx\bC\bP_1$.
All left invariant polarizations but one can be obtained by the following construction.
The (left) action of $G$ on $\bR$ extends to an affine action on $\bC$.
Fix $s\in\bC$, let
$$
f_s(g)=gs\qquad\text{ and }\qquad Q(s)=(f_{s*})^{-1} T^{1,0}\bC.\tag1
$$
E.g., $Q(i)$ gives the usual complex structure on $G\cong\bR^2$, while for $s$ real $Q(s)$ is a real polarization whose leaves are straight lines of slope $-1/s$.
Equivalently, with $H_s=\{g\in G\colon gs=s\}$, $s\in\bR$, the leaves of $Q(s)$ are left translates of $H_s$.
We write $Q(\infty)$ for the exceptional, real polarization,
whose leaves have slope $0$.
Since $f_{gs}=f_s\circ R_g$, (1) implies that $R_g\colon (G,Q(gs))\to (G,Q(s))$ is a polarized map.

Define the character $\chi\colon G\to\bR$ by $\chi(g)=b$ if 
$gt=a+bt$, and let $G^\rho=\{g\in G\colon |\chi(g)|\leq\rho\}$, 
$0\le\rho\leq\infty$.
Thus $G^\rho\subset G$ is a normal sub--semigroup if 
$\rho\leq 1$.
Everything that we discussed in this section applies to $G^1$ instead of $G$ as well.

\subhead 4.\ Riemannian manifolds\endsubhead
Let $M$ be a complete Riemannian manifold, $\dim M=m>0$, 
and $N$ the manifold of its geodesics.
The map $N\ni x\mapsto\dot x(0)\in TM$ identifies $N$ and $TM$.
We call this the canonical identification.
Following Souriau's philosophy in [So], we shall mostly avoid using this identification, though; for any $t_0\in\bR$ the map $x\mapsto\dot x(t_0)$ would define just as natural an identification anyway.
There is still a canonical symplectic form $\omega$ on $N$, see
[W, 2.3]. Under any of the above identifications it corresponds 
to the form $\sum dq_i\wedge dp_i$ (written in terms of
the usual local coordinates on $TM\approx T^*M$.) Note the sign 
convention used for $\omega$.

Composition of geodesics with affine transformations 
$\bR\to\bR$ defines a right action of $G$ on $N$. We shall consider adapted 
polarizations on $G^1$--invariant open subsets $V\subset N$.
The domains on which the adapted complex 
structures of [GS, LSz1] were defined correspond to the set of geodesics 
of speed $<r$, with $r\in (0,\infty]$ (but since then adapted complex 
structures on more general invariant sets have turned out to be 
of importance, see [FHW] and references there).
Put $\Omega_x(g)=A_g(x)=x\circ g$ for $x\in N, g\in G$. Then
$$
A^*_g\omega=\chi(g)\omega.\tag2
$$
When $gt=a+t$, (2) expresses the fact that the geodesic flow 
preserves $\omega$, and when $gt=bt$, (2) holds because in
the canonical identification $N\approx TM$ in local coordinates
$A_g$ is given by $(q_j,p_j)\mapsto (q_j,bp_j)$. Since $G$ is
generated by translations and dilations, (2) holds for all $g\in G$.

By Definition 1, given a polarization $Q$ of $G^1$, a 
polarization $P$ of $V$ is adapted to $Q$ if 
$\Omega_x\colon (G^1, Q)\to (V, P)$ is polarized for every 
$x\in V$.
If $Q$ is left invariant and $0<\rho<\infty$, this is equivalent 
to saying that 
$\Omega_x\colon (G^\rho,Q)\to (V,P)$ is polarized for every 
$x\in VG^{1/\rho}$.

\proclaim{Theorem 2}(a) If a nonempty, $G^1$--invariant open
$V\subset N$ has a polarization $P$ adapted to a polarization 
$Q$ of $G^1$, then $Q$ is left invariant.
If $Q=Q(s)$, $s\in\bC$, then it determines $P$ uniquely.\newline
\phantom{Th}(b) If $M$ is a closed analytic Riemannian 
manifold, then there 
is a $G^1$--invariant open $V\subset N$, containg all zero speed
geodesics, such that $VG^{1/|\text{Im}\,s|}$ has a polarization 
$P$ adapted to $(G^1,Q(s))$, for every $s\in\bC$. The same is
true if $M$ is not necessarily closed, but modulo its 
isometry group it is compact.
\endproclaim

It is straightforward that, after the canonical identification 
$N\cong TM$, a polarization adapted to $Q(i)$ becomes the adapted
complex structure of \cite{LSz1, Definition 4.1}.

\demo{Proof}Suppose $(V,P)$ is adapted to $Q$.
Fix a nonconstant geodesic $x\in V$; then 
$\Omega_x\colon (G^1,Q)\to (V,P)$ is a polarized immersion.
Since $\Omega_{x\circ g}=\Omega_x\circ L_g$ is also polarized for $g\in G^1$, it follows that $L_g\colon (G^1,Q)\to (G^1,Q)$ is polarized, i.e., $Q$ is left invariant.

Let now $Q=Q(s),\ s\in\bC$.
Uniqueness and existence of the adapted polarization for $s=i$ is the content of \cite{LSz1, Theorem 4.2 and Sz1, Theorem 2.2}.
To go further, with $g\in G$ and $x\in V$ consider the 
commutative diagram
$$
\CD
VG^{1/|\chi(g)|} @>A_g>>  V\\
@A\Omega_x AA @A\Omega_x AA\\
G^{1/|\chi(g)|} @>R_g>> G^1,
\endCD
$$
and recall that $R_g\colon (G, Q(gs))\to (G,Q(s))$ is polarized.
Now $A_g\colon VG^{1/|\chi(g)|}\to V$ is a diffeomorphism if 
$\chi(g)\neq 0$, and the diagram implies that $A_g$ pulls back any 
$Q(gs)$--adapted polarization to a $Q(s)$--adapted polarization.
Therefore uniqueness and existence for $s=i$ implies the same for any $s\in\bC\backslash\bR$.
Concretely, if $gs=i$, then $\chi(g)=1/\text{Im}\,s$, and so 
if $V$ admits a $Q(i)$--adapted polarization, then 
$VG^{1/|\text{Im}\,s|}$ will admit a $Q(s)$--adapted polarization.

Finally, let $s\in\bR$ and $\pi_s\colon N\to M$ be given by $\pi_s(x)=x(s)$.
As said, the leaves of $Q(s)$ are left translates of the sub--semigroup $H_s\subset G$.
That $P$ is adapted to $Q(s)$ means it is tangential to the orbits $\Omega_x(H_s)$.
As $x$ ranges over a fiber $\pi_s^{-1}y$, these orbits all pass through the constant geodesic $\equiv y$, that we denote $\overline y$.
Their tangents at $\overline y$ form the vertical tangent space $T_{\overline y}(\pi_s^{-1} y)$, which then must agree with $P_{\overline y}$.
Furthermore, the vector field generating the $H_s$--action, being tangent to the orbits, is a section of $P$.
Since $P$ is involutive, it must be invariant under the action of $H_s\cap G^1$.
But $\pi_s$ is also invariant, hence for any $x\in V$ and $\overline y$ as above, $\pi_{s*} P_x=\pi_{s*} P_{\overline y}=0$.
Therefore $P$ consists of the tangent spaces to the fibers of $\pi_s$, and is unique.
It is straightforward that, conversely, the tangent spaces to the fibers form a polarization of $N$, adapted to $(G^1,Q(s))$ (and to $(G,Q(s))$).
\enddemo

The proof also gave the following

\proclaim{Corollary 3}Let $s\in\bC$ and $g\in G$.
If a $G^1$--invariant $V\subset N$ admits a $Q(s)$--adapted
polarization $P(s)$, then $VG^{1/|\chi(g)|}$ admits a 
$Q(gs)$--adapted polarization $P(gs)$, and
$$
A_g\colon (VG^{1/|\chi(g)|}, P(gs))\to (V, P(s))
$$
is polarized, in fact a polarized isomorphism when $\chi(g)\neq 0$.
\endproclaim

We shall continue using $P(s)$ for the $Q(s)$--adapted polarization 
on $V\subset N$, whenever it exists.
Let $v(x)$ denote the speed squared of a geodesic $x\in N$.

\proclaim{Theorem 4}For $s\in\bC\backslash\bR$ let $\partial_s,\overline\partial_s$ denote the complex exterior derivations 
for the complex structure $P(s)$ on $V$ (if this latter exists).
The symplectic form 
$\omega$ on $V\subset N\cong TM$ is given by
$$
i\omega=(\roman{Im}\,s)\overline\partial_s\partial_s v.
$$
In particular, $\omega$ is a positive or negative $(1,1)$--form depending on whether $\roman{Im}\,s>0$ or $<0$.
\endproclaim

\demo{Proof}When $s=i$, the claims are in \cite{LSz1, Corollary 5.5 and Theorem 5.6}. (Note that $E$ and $\Omega$ there correspond to
our $v/2$, resp. $-\omega$.)
Hence the general case follows by Corollary 3, because 
$A_g^*v=\chi(g)^2 v$ and by (2), $A_g^*\omega=\chi(g)\omega$.
\enddemo

\subhead 5.\ The canonical bundle\endsubhead
In the set up of Section 4, let $s\in\Bbb C\setminus\Bbb R$
and consider a $G^1$--invariant open
$V\subset N$ that admits a $Q(s)$--adapted polarization $P(s)$,
a complex structure. Its canonical bundle, the holomorphic 
line bundle $K\to V$ of $(m,0)$--forms, has a Hermitian metric
$h^K$ defined by
$$
h^K(\theta)\omega^m(x)=i^{m^2}m!\theta\wedge\bar\theta,\qquad
\theta\in K_x,\quad x\in V.\tag3
$$
In this subsection we compute $h^K$, 
something that is needed for purposes of
quantization. We start by recalling certain constructions and
results from [LSz1, Sz2]. Denote by id$\in G$ the identity 
transformation $\Bbb R\to\Bbb R$.
 
The action of $G$ on $N$ induces an action on $TN$ and
$\Bbb C\otimes TN$, denoted
$(\xi,g)\mapsto \xi g$. Let $x\in V$. Any 
$\xi\in T_xV$ can be decomposed into $(1,0)$ and $(0,1)$ 
components with respect to the structure $P(s)$:
$\xi=\xi^{1,0}+\xi^{0,1}$. If $J:TV\to TV$ denotes the complex 
structure operator for $P(s)$, then $\xi^{1,0}=(\xi-iJ\xi)/2$.
The map $g\mapsto(\xi g)^{1,0}$ is holomorphic as a map
$(G^1,Q(s))\to T^{1,0}(V,P(s))$ (in the sense that it has a 
holomorphic extension to a neighborhood of $G^1$),
see [LSz1, Proposition 5.1]. 

Now consider two $m$--tuples 
$\xi_1,\ldots$ and $\eta_1,\ldots\in T_xN$, and assume that
the $\xi^{1,0}_j$ are linearly independent.
Those $g$ for which $(\xi_jg)^{1,0}$
 are linearly dependent form 
a discrete subset $\Delta\subset G$. For 
$g\in G^0\setminus\Delta$ the 
$\xi_jg$ are also independent. Since when $g\in G^0$, the
vectors $\xi_jg,\eta_jg$ are tangential to the $m$--dimensional
manifold of zero speed geodesics, on $G^0\setminus\Delta$
there is a smooth real $m\times m$--matrix valued function
$\phi^0=(\phi^0_{jk})$ such that
$$
\eta_jg=\sum_k\phi^0_{jk}(g)\xi_kg,\qquad g\in 
G^0\setminus\Delta.\tag4
$$
Further, there is a meromorphic $m\times m$--matrix valued
function $\phi=(\phi_{jk})$ on $(G^1,Q(s))$, 
with poles restricted to $\Delta$, such that
$$
(\eta_jg)^{1,0}=\sum_k\phi_{jk}(g)(\xi_kg)^{1,0}.\tag5
$$
By (4) and (5), $\phi$ is the analytic continuation of $\phi^0$.

\proclaim{Theorem 5} Suppose $x\in V$ and
$\xi_1,\ldots,\eta_m\in T_xV$ form a symplectic basis:
$$
\omega(\xi_j,\xi_k)=
\omega(\eta_j,\eta_k)=0,\quad
\omega(\xi_j,\eta_k)=\delta_{jk},\qquad 1\le j,k\le m.
$$
If $\phi$ is as in (5), then for $\theta\in K_x$
$$ 
h^K(\theta)=
2^m|\theta(\xi_1,\ldots,\xi_m)|^2
\det\roman{Im}\,\phi(\roman{id}).
\tag6
$$
\endproclaim
\demo{Proof} With $\zeta_j=
\sum_k\text{Im}\,\phi_{jk}(\text{id})\xi_k$
$$
\multline
\theta\wedge\bar\theta(\zeta_1^{1,0},\ldots,\zeta_m^{1,0},
\eta_1^{0,1},\ldots,\eta_m^{0,1})=
\theta(\zeta_1^{1,0},\ldots)\overline{\theta(\eta_1^{1,0},\ldots)}\\
=\theta(\zeta_1,\ldots)\overline{\theta(\eta_1^{1,0},\ldots)}
=|\theta(\xi_1,\ldots)|^2\det\text{Im}\,\phi(\text{id})
\det\bar\phi(\text{id}).
\endmultline\tag7
$$
Taking real parts in (5) gives 
$\eta_j=\sum_k\text{Re }\phi_{jk}(\text{id})\xi_k+J\zeta_j$. Thus
$$\align
2\omega(\zeta^{1,0}_j,\eta^{0,1}_l)=2\omega(\zeta^{1,0}_j,\eta_l)
=\omega(\zeta_j-&iJ\zeta_j,\eta_l) \\
=\omega\big(\sum_k\text{Im}\,\phi_{jk}(\text{id})
\xi_k,\eta_l\big)-
i&\omega\big(\eta_j-\sum_k\text{Re }\phi_{jk}(\text{id})
\xi_k,\eta_l\big) \\
=\text{Im}\,\phi_{jl}(\text{id})+i&\text{Re }\phi_{jl}(\text{id})=i\bar\phi_{jl}(\text{id}),
\endalign
$$
and
$$\align
\omega^m(\zeta_1^{1,0},\ldots,\zeta_m^{1,0},
\eta_1^{0,1},\ldots,\eta_m^{0,1})&=
m!i^{m(m-1)}\det\bigl(\omega(\zeta^{1,0}_j,\eta^{0,1}_l)\bigr) 
\\ &=m!i^{m^2}2^{-m}\det\bar\phi(\text{id}).
\endalign$$
Comparing this with (3) and (7) yields (6).\enddemo

\subhead 6.\, The family of adapted polarizations\endsubhead
Finally we shall construct a polarized fibration $Z\to\bC$ whose fibers represent the various $(V,P(s))$.
With $s\in\bC$, $x\in N$ consider the embeddings
$$
i^x\colon\bC\ni s\mapsto (s,x)\in
\bC\times N,\qquad j^s\colon N\ni x\mapsto (s,x)\in\bC\times N.
\tag8
$$
Also, let $\pi\colon\bC\times N\to\bC$ denote the projection.

\proclaim{Theorem 6}Suppose that a $G^1$--invariant open
$V\subset N$ admits the adapted polarization (complex structure) $P(i)$. \newline
\phantom{Th}(a) On 
$Z=\{(s,x)\in\bC\times N\colon x\in VG^{1/|\roman{Im}\,s|}\}$
there is a unique polarization $P$ such that the maps
$$
i^x\colon ((i^x)^{-1} Z, T^{1,0}\bC)\to (Z,P),\qquad  
j^s\colon (VG^{1/|\roman{Im}\,s|},P(s))\to (Z,P)
$$
are polarized for all $s\in\bC$, $x\in N$. With this 
$P$, $\pi\colon (Z,P)\to (\bC,T^{1,0}\bC)$ is polarized,
and 
$(Z\setminus\pi^{-1}\Bbb R, P)=Z_0$ is a complex manifold.
\newline
\phantom{Th}(b)
Let $\partial$, $\bar\partial$ denote the complex exterior
derivations on $Z_0$, and $\tilde\omega$ 
the pullback of $\omega$ along the map $(s,x)\to x$. Then
$$
i\tilde\omega=\bar\partial\partial(v\roman{Im}\,s)\qquad
\text{on }Z_0,\tag9
$$
$v\roman{Im}\,s$ is plurisub/superharmonic if 
$\roman{Im}\,s>0$, resp. $<0$, and satisfies the Monge--Amp\`ere
equation $\roman{rk}\,\bar\partial\partial(v\roman{Im}\, s)=m$. 
\newline
\phantom{Th}(c)
Finally, endow $(\bC, T^{1,0}\bC)\times (V,P(i))$ with the 
product complex structure.
Then the map $\Phi\colon Z\to\bC\times V$ given by
$$
\Phi(s,x)=(s,x\circ g),\qquad\text{where}\quad gi=f_i(g)=s,\tag10
$$
is polarized, and in fact restricts to a biholomorphism 
$Z_0\to(\bC\backslash\bR)\times V$.
\endproclaim

\demo{Proof}Since the range of $i^x_*$ and $j^s_*$ together span 
$T(\bC\times N)$, the polarization $P$ in question is unique, and must be given by
$$
P_{(s,x)}=i^x_* T^{1,0}_s\bC\oplus j^s_* P(s)_x,\quad (s,x)\in Z.\tag11
$$
In view of Corollary 3 this formula defines a subbundle $P\subset\bC\otimes TZ$.
Our $P$ has rank $m+1$ all right, but is it involutive? To decide,
first we check that $\Phi$ in (10) is polarized, i.e.~$\Phi_*$ maps $P$ into $T^{1,0}(\bC\times V)$.
With notation introduced earlier
$$
(\Phi\circ i^x)(s)=(s,(\Omega_x\circ f_i^{-1})(s)),
\qquad (\Phi\circ j^s)(x)=(s,A_g x).
$$
Now $\Omega_x\colon (G^{r/\sqrt{v(x)}}, Q(i))\to (V, P(i))$ 
is holomorphic by the definition of $P(i)$ and by the observation 
preceding Theorem 2; also 
$f_i\colon (G,Q(i))\to\bC$ is biholomorphic by the definition of 
$Q(i)$. Therefore $\Phi\circ i^x$ is holomorphic. Similarly 
$\Phi\circ j^s\colon (VG^{1/|\text{Im}\,s|},P(s))\to
(\bC\times V, T^{1,0}(\bC\times V))$ is polarized by Corollary 3 
($s$ there corresponds to $i$ here, though).
Putting these and (11) together, we see 
$\Phi_* P\subset T^{1,0} (\bC\times V)$ indeed.

Since $T^{1,0}(\bC\times V)$ is involutive and $\Phi$ is a diffeomorphism over $Z_0$, it follows that $P$ is involutive over 
$Z_0$, which therefore is a complex manifold. By density, 
$P$ is involutive over all of $Z$.
That $\pi$ is polarized is obvious, so (a) and (c) have been proved.
As to (9), the two sides restricted to the fibers of $\pi$ agree 
by Theorem 4. Tangents
to the fibers of $(s,x)\to x$, the ``horizontal fibers'', constitute
the kernel of of $i\tilde\omega$; to finish the proof it will
suffice to show the same for 
$\text{Ker }\bar\partial\partial(v\text{Im}\,s)$. The restriction
of $\bar\partial\partial(v\text{Im}\,s)$ to the horizontal fibers 
is certainly $0$, since $v$ restricts to a constant and $i^x$
is holomorphic; but that is not quite enough. It will be necessary 
to compute $\bar\partial\partial(v\text{Im}\,s)$,
that we do by pulling it back along $\Phi^{-1}$.

If in (10) $gt=a+bt$, then $b=\text{Im}\,s$. Hence 
$(\Phi^{-1})^*(v\text{Im}\,s)=v/\text{Im}\,s$. (With a slight abuse 
of notation, $v$ stands for both a function on $N$ and its
pull back along the projection $\Bbb C\times N\to N$. Also,
$\text{Im}\,s$ stands for a function on $\Bbb C\times N$.)
On $(\Bbb C\setminus\Bbb R)\times V$
$$ 
\bar\partial\partial\frac{v}{\text{Im}\,s}=
\frac{\bar\partial\partial v}{\text{Im}\,s}+
\frac{d\bar s\wedge\partial v-\bar\partial v\wedge ds}
{2i(\text{Im}\,s)^2}+\frac{vd\bar s\wedge ds}{2(\text{Im}\,s)^3}.
\tag12
$$
In computing $\partial v$, $\bar\partial v$, the operators
corresponding to the complex structure $P(i)$ are to be used. 
From (12) $(\bar\partial\partial v/\text{Im}\,s)^{m+1}$ equals 
$$
\multline
(m+1)\big(\frac{\bar\partial\partial v}{\roman{Im}\,s}\big)^m\wedge
\frac{v\,d\bar s\wedge ds}{2(\roman{Im}\,s)^3}-\binom {m+1}2
\big(\frac{\bar\partial\partial v}{\roman{Im}\,s}\big)^{m-1}\wedge
\frac{(d\bar s\wedge\partial v-\bar\partial v\wedge ds)^2}
{4(\roman{Im}\,s)^4}\\
=\frac{(m+1)d\bar s\wedge ds}{4(\roman{Im}\, s)^{m+3}}\wedge
\big(2v(\bar\partial\partial v)^m-
m(\bar\partial\partial v)^{m-1}\wedge\bar\partial v\wedge
\partial v\big).
\endmultline
$$
By [LSz1, (5.10)], where $E=v/2$, the last expression
vanishes. As $-i\bar\partial\partial (v/\roman{Im}\, s)$
is definite along the fibers of $\pi$, its signature
is $m$ pluses (or minuses) and one $0$, and the same
holds for $-i\bar\partial\partial (v\roman{Im}\, s)$ on $Z_0$.
In particular, $(v\roman{Im}\, s)$ is plurisub/superhar\-mon\-ic,
and because $\bar\partial\partial(v\roman{Im}\, s)$ vanishes 
on the horizontal fibers, its kernel consists of the tangents
to the horizontal fibers. This then proves (9) and the rest of (b).

\enddemo

\enddocument

\Refs
\widestnumber\key{FMMN1}
\ref\key{B}\by R.~Bielawski\paper Complexification and hypercomplexification of manifolds with a linear connection
\jour Internat.~J.~Math.\vol 14\yr 2003\pages 813--824\endref

\ref\key{FHW}\manyby G. Fels, A. Huckleberry, J. Wolf\book
Cycle spaces of flag domains. A 
complex geometric viewpoint. \publ Birkh\"auser \publaddr
Boston, MA \yr 2006\endref

\ref\key{FMMN1}\manyby C.~Florentino, P.~Matias, J.~Mour\~ao, J.P.~Nunes\paper Geometric quantization, complex structures and the coherent state transform\jour J.~Funct.~Anal.\vol 221\yr 2005\pages 303--322\endref

\ref\key{FMMN2}\bysame \paper On the BKS pairing for K\"ahler quantizations of the cotangent bundle of a Lie group
\jour J.~Funct.~Anal.\vol 234\yr 2006\pages 180--198\endref

\ref\key{GS}\by V.~Guillemin, M.B.~Stenzel\paper Grauert tubes and the homogeneous Monge--Amp\`ere equation\jour J.~Diff.~Geom\vol 34\yr 1991\pages 561--570\endref

\ref\key{H1}\manyby B.C.~Hall\paper The Segal--Bargmann ``coherent state'' transform for compact Lie groups\jour J.~Funct.~Anal.\vol 122\yr 1994\pages 103--151\endref

\ref\key{H2}\bysame \paper Geometric quantization and the generalized Segal--Bargmann transform for Lie groups of compact type\jour Comm.~Math.~Phys.\vol 226\yr 2002\pages 233--268\endref

\ref\key{HK}\by B.C.~Hall, W.D.~Kirwin\paper Adapted complex structures and the geodesic flow, arxiv: 0811.3083\endref

\ref\key{L}\by L.~Lempert\paper Complex structures on the tangent bundle of Riemannian manifolds\inbook Complex Analysis and Geometry\publ Plenum\publaddr New York, N.Y.\yr 1993\pages 235--251\endref

\ref\key{LSz1}\manyby L.~Lempert, R.~Sz\H{o}ke\paper Global solutions of the homogeneous complex Monge--Amp\`ere equation and complex structures on the tangent bundle of Riemannian manifolds\jour Math.~Ann.\vol 290\yr 1991\pages 689--712\endref

\ref\key{LSz2}\bysame\paper Uniqueness in geometric quantization, manuscript\endref


\ref\key{So}\by J.-M.~Souriau\book Structure of dynamical systems.
A symplectic view of physics\publ Birkh\"auser\publaddr Boston, MA\yr 1997\endref

\ref\key{Sz1}\manyby R.~Sz\H{o}ke\paper Complex structures on tangent bundles of Riemannian manifolds\jour Math.~Ann.\vol 291\yr 1991\pages 409--428\endref

\ref\key{Sz2}\bysame\paper Automorphisms of certain Stein manifolds
\jour Math. Z.\vol 219\yr 1995\pages 357--385\endref

\ref\key{Sz3}\bysame\paper Canonical complex structures associated to connections and complexifications of Lie groups\jour Math.~Ann.\vol 328\yr 2004\pages 553--591\endref

\ref\key{W}\by N.M.J.~Woodhouse\book Geometric Quantization
\bookinfo 2nd edition\publ Clarendon Press\publaddr Oxford, UK\yr 1992\endref

\endRefs
\bye